\documentclass[aps,pre,twocolumn,groupedaddress,amsmath,amssymb]{revtex4-2}

\usepackage{graphicx}
\usepackage{bm}
\usepackage{color}

\newcommand{\ctext}[1]{\raise0.2ex\hbox{\textcircled{\scriptsize{#1}}}}

\begin{document}


\title{Statistics for stochastic differential equations and approximations of resolvent}


\author{Jun Ohkubo}
\email[]{johkubo@mail.saitama-u.ac.jp}
\affiliation{
Graduate School of Science and Enginnering, Saitama University,
255 Shimo-Okubo, Sakura-ku, Saitama, 338-8570, Japan
}




\begin{abstract}
The numerical evaluation of statistics plays a crucial role in statistical physics and its applied fields. It is possible to evaluate the statistics for a stochastic differential equation with Gaussian white noise via the corresponding backward Kolmogorov equation. The important notice is that there is no need to obtain the solution of the backward Kolmogorov equation on the whole domain; it is enough to evaluate a value of the solution at a certain point that corresponds to the initial coordinate for the stochastic differential equation. For this aim, an algorithm based on combinatorics has recently been developed. In this paper, we discuss a higher-order approximation of resolvent, and an algorithm based on a second-order approximation is proposed. The proposed algorithm shows a second-order convergence. Furthermore, the convergence property of the naive algorithms naturally leads to extrapolation methods; they work well to calculate a more accurate value with fewer computational costs. The proposed method is demonstrated with the Ornstein-Uhlenbeck process and the noisy van der Pol system.
\end{abstract}


\maketitle

\section{Introduction}

Numerical time integration is a common theme in various research areas of computer simulation. When considering a system with noise, one usually focuses on statistics such as mean and variance. These statistics play a crucial role in statistical physics and various applications such as data assimilation and filtering; for example, see Ref.~\cite{Evensen_book}. It is common to use Monte Carlo simulations for evaluating the statistics numerically. For example, the Euler-Maruyama approximation is well known for numerical simulations of stochastic differential equations~\cite{Kloeden_Platen_book}. Stochastic differential equations have many applications ranging from time series data analysis to control problems. Furthermore, rapid evaluation of statistics is necessary for real-time processing. However, generations of random numbers and samplings for accuracy require high computational costs.

A stochastic differential equation with multiplicative Gaussian white noise has a corresponding partial differential equation, i.e., a Fokker-Planck equation \cite{Gardiner_book, Risken_book}. In addition, as explained later, a different type of partial differential equation, the backward Kolmogorov equation, is also available to discuss the statistics for the stochastic differential equation. Since the backward Kolmogorov equation is a partial differential equation, it is sometimes hard to obtain its numerical solution. There are several approaches; for example, the conventional space-time discretization methods are available \cite{Knabner_and_Angermann_book}. Recent studies with low-rank approximations give new numerical methods for partial differential equations; for example, see Refs.~\cite{Dolgov2012,Rodgers2020,Einkemmer2021,Dektor2021,Dektor2021b}. Some works employ deep neural networks to obtain solutions to partial differential equations; for example, see Refs.~\cite{E2017,Raissi2019,E2022}. The other approach is moment-closure methods; there are many references for this topic, ranging from stock price models \cite{Singer2006} to population models \cite{Ekanayake2010}; see the review in Ref.~\cite{Kuehn2016}. However, this moment-closure approach needs a kind of truncation.

Spatial-discretization methods, Galerkin methods, and collocation methods are famous as numerical approaches for partial differential equations \cite{Knabner_and_Angermann_book}. These approaches lead to coupled ordinary differential equations to approximate an original partial differential equation. Then, we employ numerical integration methods such as the explicit Euler method, the Crank-Nicolson method, and the Runge-Kutta method. Note that these numerical approaches aim to obtain the solution of the partial differential equation on the whole domain.

The focus here is the evaluation of statistics for stochastic differential equations. However, there is no need to evaluate the complete solution of the backward Kolmogorov equation. It is sufficient to obtain only the value of the solution at the point corresponding to the initial condition of the stochastic differential equations; we will see it later. Hence, there is no need to use conventional numerical methods for partial differential equations. For this purpose, Ref.~\cite{Ohkubo2021} gave an efficient algorithm to evaluate the statistics for short-time ranges with the aid of resolvent and combinatorics. Algorithms based on combinatorics are powerful, and some recent works in statistical physics focus on them. For example, there are discussions based on combinatorics to compute the Mori-Zwanzig memory integral in generalized Langevin equations \cite{Amati2019,Zhu2020}. However, as pointed out in Ref.~\cite{Ohkubo2021}, one of the problems of the algorithm based on combinatorics is the slow convergence speed in the iterative procedure. Is it possible to achieve more rapid convergence? Furthermore, the evaluation is also important in revealing the nonlinear dynamics of the system. Recently, methods based on the Koopman operator have made significant progress. For example, the spectral properties of the Koopman operator have been compared to the nonlinear dynamics in physical systems~\cite{Mezic2004}; the method was extended to deal with the nonlinearities more adequately~\cite{Williams2015}. These methods are data-driven, and one can recover even system equations from time series data~\cite{Klus2020}. However, it would be beneficial to develop not only data-driven approaches to the Koopman operators but also equation-based ones to understand the nonlinear behavior of the system from a physical perspective. Indeed, Ref.~\cite{Ohkubo2022} clarified the relationship between the backward Kolmogorov equation and the Koopman operator; the approach based on combinatorics provided a practical numerical tool to evaluate the matrix elements of the Koopman operator. In this context, numerical algorithms for faster and more efficient evaluation of statistics in stochastic processes are broadly beneficial for studies of nonlinear dynamics.

The present paper discusses algorithms based on combinatorics, especially higher-order approximations of the resolvent, which gives us a more efficient numerical method. Algorithms based on combinatorics sometimes use the inverse of a matrix of infinite dimension. Hence, their local approximation is employed. Here, the local approximation means neglecting higher-order `events' to evaluate the matrix elements. The 2nd-order local approximation of resolvent is explicitly derived, and the proposed algorithm shows 2nd-order convergence. By contrast, the previous algorithm has 1st-order convergence. Thus, the proposed algorithm can compute statistics faster and more efficiently. Furthermore, the convergence property of the proposed algorithm naturally leads to an extrapolation method, which reduces the computational cost. We demonstrate the proposed algorithm by applying it to the Ornstein-Uhlenbeck process and the noisy van der Pol system.

The construction of the present paper is as follows. In Sec.~II, we briefly review algorithms based on combinatorics. We will see that a specific type of basis expansion avoids the need to solve the backward Kolmogorov equation on the whole domain. Section III gives the first main contribution of the present paper; we derive a 2nd-order approximation of the matrix inversion derived from the time-evolution operator. We also yield numerical demonstrations of the proposed method for Ornstein-Uhlenbeck processes and noisy van der Pol systems. In Sec.~IV, the 2nd-order convergence property leads to an extrapolation method, which allows us to evaluate the statistics with less computational cost. Section V gives some concluding remarks.

\section{Combinatorics for Evaluating Statistics}
\label{sec_combinatorics}

\subsection{Stochastic differential equation and Fokker-Planck equation}
\label{subsec_sde}

We briefly review the basics of a stochastic differential equation and the corresponding Fokker-Planck equation. For details, see, for example, Ref.~\cite{Gardiner_book} and Ref.~\cite{Risken_book}.

Consider a $D$-dimensional state space $\mathcal{M} \subseteq \mathbb{R}^D$ and a vector of stochastic variables $\bm{X} \in \mathcal{M}$. We consider continuous time-evolution of the state vector $\bm{X}$, and it obeys the following stochastic differential equation:
\begin{align}
d\bm{X} = \bm{a}(\bm{X}) d t + \bm{B}(\bm{X}) d\bm{W}(t),
\label{eq_SDE}
\end{align}
where $\bm{a}(\bm{X})$ is a vector of drift coefficient functions, and $\bm{B}(\bm{X})$ is a matrix of diffusion coefficient functions. Although the coefficient functions can be time-dependent, we here discuss only time-independent cases for simplicity. The components of the vector of Wiener processes $\bm{W}(t)$ satisfy
\begin{align}
&\mathbb{E}\left[ W_i(t) \right] = 0, \\
&\mathbb{E}\left[ \Big(W_i(t) - W_i(s) \Big)  \Big(W_j(t) - W_j(s) \Big)\right] = (t-s)\delta_{ij},
\end{align}
for $0 \le s \le t$. We specify the initial condition for the stochastic differential equation as $\bm{X}(0) = \bm{x}_\mathrm{ini}$. The stochastic differential equation in \eqref{eq_SDE} has a corresponding Fokker-Planck equation, which describes the time-evolution of the probability density function $p(\bm{x},t)$:
\begin{align}
\frac{\partial}{\partial t} p(\bm{x},t) = \mathcal{L} \, p(\bm{x},t),
\label{eq_Fokker_Planck}
\end{align}
where
\begin{align}
\mathcal{L}
= - \sum_{i} \frac{\partial}{\partial x_i} a_i(\bm{x}) 
+ \frac{1}{2} \sum_{i,j} \frac{\partial^2}{\partial x_i \partial x_j} 
\left[ \bm{B}(\bm{x}) \bm{B}(\bm{x})^\mathrm{T} \right]_{ij}
\end{align}
is the time-evolution operator for the Fokker-Planck equation. The initial condition at the initial time $t = 0$ is 
\begin{align}
p(\bm{x}, t=0) = \delta(\bm{x} - \bm{x}_\mathrm{ini}),
\end{align}
where $\delta(\cdot)$ is the Dirac delta function.

By solving the Fokker-Planck equation, we obtain a probability density function or a transition probability for $\bm{X}$ at time $t$. For example, provided the initial condition is $\bm{x}_\mathrm{ini}$, the $\bm{\alpha}$-th order moment at time $t = T$ is evaluated as
\begin{align}
\mathbb{E}\left[  \bm{X}(T)^{\bm{\alpha}} | \bm{X}(0) = \bm{x}_\mathrm{ini}  \right] 
= \int_{\mathcal{M}} \bm{x}^{\bm{\alpha}} p(\bm{x},T) d\bm{x},
\label{eq_expectation_basic}
\end{align}
where $\alpha_d \in \mathbb{N}_0$ for $d = 1, \dots, D$ and $\bm{x}^{\bm{\alpha}} = x_1^{\alpha_1} \cdots x_D^{\alpha_D}$.

\subsection{Backward Kolmogorov equation}
\label{subsec_backward_Kolmogorov}

When one wants to focus on the statistics for the stochastic differential equation, the corresponding backward Kolmogorov equation is available. The following transformation for \eqref{eq_expectation_basic} leads to the corresponding backward Kolmogorov equation formally; repeating the integration by parts, we have
\begin{align}
\int_{\mathcal{M}} \bm{x}^{\bm{\alpha}} p(\bm{x},T) d\bm{x} 
&= \int_{\mathcal{M}}
\bm{x}^{\bm{\alpha}}
\left( e^{\mathcal{L} T}  \delta(\bm{x} - \bm{x}_\mathrm{ini}) \right) 
d\bm{x}\nonumber \\
&= \int_{\mathcal{M}}
\left( e^{\mathcal{L}^\dagger T}  \bm{x}^{\bm{\alpha}} \right)
\delta(\bm{x} - \bm{x}_\mathrm{ini}) 
d\bm{x} \nonumber \\
&= \int_{\mathcal{M}} \varphi_{\bm{\alpha}}(\bm{x},T) \delta(\bm{x} - \bm{x}_\mathrm{ini}) d\bm{x} \nonumber \\
&= \varphi_{\bm{\alpha}}(\bm{x}_\mathrm{ini},T),
\label{eq_expectation_basic2}
\end{align}
where
\begin{align}
\mathcal{L}^\dagger
= \sum_{i} a_i(\bm{x}) \frac{\partial}{\partial x_i} 
+ \frac{1}{2} \sum_{i,j} \left[ \bm{B}(\bm{x}) \bm{B}(\bm{x})^\mathrm{T} \right]_{ij} \frac{\partial^2}{\partial x_i \partial x_j},
\label{eq_adjoint_operator}
\end{align}
and $\varphi_{\bm{\alpha}}(\bm{x},t)$ is the solution of the following equation:
\begin{align}
\frac{d}{dt} \varphi_{\bm{\alpha}}(\bm{x},t) = \mathcal{L}^\dagger \varphi_{\bm{\alpha}}(\bm{x},t),
\label{eq_time_evolution_adjoint}
\end{align}
which leads to $\varphi_{\bm{\alpha}}(\bm{x},T) = \exp(\mathcal{L}^\dagger T) \varphi_{\bm{\alpha}}(\bm{x},0)$ in \eqref{eq_expectation_basic2}. Note that the initial condition for \eqref{eq_time_evolution_adjoint} should be
\begin{align}
\varphi_{\bm{\alpha}}(\bm{x},0) = \bm{x}^{\bm{\alpha}}.
\label{eq_initial_for_naive_backward_Kolmogorov}
\end{align}
The partial differential equation \eqref{eq_time_evolution_adjoint} is the backward Kolmogorov equation for \eqref{eq_SDE}.

There are some comments. First, \eqref{eq_expectation_basic2} indicates that there is no need to evaluate the solution $\varphi_{\bm{\alpha}}(\bm{x},T)$ on the whole domain. Only a value of $\varphi_{\bm{\alpha}}(\bm{x}_\mathrm{ini},T)$ is enough to calculate the value of the target statistic for the initial condition $\bm{x}_\mathrm{ini}$. Second, when the time-evolution operator in \eqref{eq_adjoint_operator} is time-dependent, we must solve the backward Kolmogorov equation \eqref{eq_time_evolution_adjoint} backwardly in time. If one employs algorithms based on combinatorics, which is the main topic of the present paper, an extension of random variables is necessary for time-dependent cases. Reference~\cite{Ohkubo2020} discussed the extension via the Ito formula. Since we consider only time-independent cases here, the following discussions become simple.

\subsection{Combinatorics}
\label{subsec_combinatorics}

As discussed below, it is beneficial to evaluate the statistics $\mathbb{E}\left[ \left(\bm{X}(T) -  \bm{x}_\mathrm{ini}\right)^{\bm{\alpha}} | \bm{X}(0) = \bm{x}_\mathrm{ini}  \right]$ instead of \eqref{eq_expectation_basic}. That is, instead of equation \eqref{eq_expectation_basic}, we consider
\begin{align}
&\mathbb{E}\left[  \left(\bm{X}(T) - \bm{x}_\mathrm{ini} \right)^{\bm{\alpha}} | \bm{X}(0) = \bm{x}_\mathrm{ini}  \right] \nonumber \\
&= \int_{\mathcal{M}} \left(\bm{x} - \bm{x}_\mathrm{ini}\right)^{\bm{\alpha}}  p(\bm{x},T) d\bm{x}.
\end{align}
Hence, it is useful to employ the following expansion:
\begin{align}
\varphi_{\bm{\alpha}}(\bm{x},t) = \sum_{\bm{n}} P_{\bm{\alpha}}(\bm{n},t) (\bm{x} - \bm{x}_\mathrm{ini})^{\bm{n}},
\label{eq_basis_expansion}
\end{align}
where $n_d \in \mathbb{N}_0$ for $d = 1, \dots, D$ and $\bm{n} = (n_1, \dots, n_D)$. $\{P_{\bm{\alpha}}(\bm{n},t)\}$ are the expansion coefficients, and $\bm{x}_\mathrm{ini}$ is the initial coordinate for the stochastic differential equation. Due to the shift of origin with $\bm{x}_\mathrm{ini}$, it is enough to evaluate only the coefficient $P_{\bm{\alpha}}(\bm{0},T)$ because
\begin{align}
\varphi_{\bm{\alpha}}(\bm{x}_\mathrm{ini},T) 
&= \sum_{\bm{n}} P_{\bm{\alpha}}(\bm{n},T) (\bm{x}_\mathrm{ini} - \bm{x}_\mathrm{ini})^{\bm{n}} \nonumber \\
&= P_{\bm{\alpha}}(\bm{0},T).
\label{eq_zero_coefficient}
\end{align}
Note that the initial condition in equation \eqref{eq_initial_for_naive_backward_Kolmogorov} changes as follows:
\begin{align}
\varphi_{\bm{\alpha}}(\bm{x},0) = \left(\bm{x} -\bm{x}_\mathrm{ini}\right)^{\bm{\alpha}}.
\end{align}

\noindent
The shift of origin reduces the computational costs largely because only a coefficient $P_{\bm{\alpha}}(\bm{0},T)$ is necessary. Of course, it is possible to recover the original statistics $\mathbb{E}\left[ \bm{X}(T)^{\bm{\alpha}} | \bm{X}(0) = \bm{x}_\mathrm{ini}  \right]$ by using the results for lower-order moments.

The basis expansion in \eqref{eq_basis_expansion} leads to coupled ordinary differential equations for $\{P_{\bm{\alpha}}(\bm{n},t)\}$ instead of the partial differential equation \eqref{eq_time_evolution_adjoint}. However, note that the time-evolution operator in \eqref{eq_adjoint_operator} leads to infinite systems for $\{P_{\bm{\alpha}}(\bm{n},t)\}$, which is intractable numerically. Instead of the direct numerical integration for the coupled ordinary differential equations, we here employ algorithms based on combinatorics for numerical evaluation as below. For example, if the first term in the r.h.s. of \eqref{eq_adjoint_operator} contains a term $x_d^2 \partial_{x_d}$, it is rewritten as
\begin{align}
x_d^2 \frac{\partial}{\partial x_d} 
= &\left(x_d - x_{\mathrm{ini},d}\right)^2 \frac{\partial}{\partial x_d} 
+ 2 \left(x_d - x_{\mathrm{ini},d}\right)  x_{\mathrm{ini},d} \frac{\partial}{\partial x_d} \nonumber \\
&+ x_{\mathrm{ini},d}^2 \frac{\partial}{\partial x_d}.
\label{eq_example_drift}
\end{align}
The action of the first term in the r.h.s. of \eqref{eq_example_drift} to the basis in \eqref{eq_basis_expansion} gives 
\begin{align}
(x_d - x_{\mathrm{ini},d})^2 \frac{\partial}{\partial x_d} (\bm{x} - \bm{x}_\mathrm{ini})^{\bm{n}} 
= n_d (\bm{x} - \bm{x}_\mathrm{ini})^{\bm{n}'},
\end{align}
where
\begin{align}
\bm{n}' = (n_1, \dots, n_d + 1, \dots,  n_D).
\end{align}
Then, there is a change in the `discrete state' with $\bm{n} \to \bm{n}'$. From the similar discussion, it is clear that the second term in the r.h.s. of \eqref{eq_example_drift} does not change the discrete state $\bm{n}$, and the third one decreases $n_d$ by one. This means that we have a discrete process for $\bm{n}$ instead of the partial differential equation with continuous variables in \eqref{eq_time_evolution_adjoint}. Note that it is possible to extend this discussion to the duality in stochastic processes \cite{Liggett_book,Giardina2009,Ohkubo2013,Jansen2014}, and numerical algorithms based on combinatorics have been developed in recent studies \cite{Ohkubo2021,Ohkubo2022}. Since there is no need to employ the discussion based on the duality in the present paper, we omit the details. 

The action of $\mathcal{L}^\dagger$ to the basis expansion in \eqref{eq_basis_expansion} is interpreted as the transition on the discrete state $\bm{n}$. Hence, it is useful to employ the matrix expression for \eqref{eq_adjoint_operator}; the matrix element for $\bm{n} \to \bm{n}'$ is written as
\begin{align}
\left[ \mathcal{L}^\dagger \right]_{\bm{n}' \bm{n}} = \sum_{r=1}^R \gamma_r(\bm{n}) \, \delta_{\bm{v}_r, \bm{n}'-\bm{n}},
\label{eq_explicit_expression_for_operetor}
\end{align}
where $R$ is the number of terms (events) in $\mathcal{L}^\dagger$, $\gamma_r(\bm{n})$ is a state-dependent coefficient for the $r$-th event, and $\bm{v}_r$ is a vector for the state change in $\bm{n}$ for the $r$-th event. The matrix expression is easily obtained from the time-evolution operator in \eqref{eq_adjoint_operator}. Then, the following simple expression leads to walks on the lattice: 
\begin{align}
e^{\mathcal{L}^\dagger T} \varphi_{\bm{\alpha}}(\bm{x},0)
&= \left(  e^{ (\mathcal{L}^\dagger T)/M} \right)^M \varphi_{\bm{\alpha}}(\bm{x},0) \nonumber \\
&\simeq \left(  \bm{1} + \frac{T}{M} \mathcal{L}^\dagger \right)^M \varphi_{\bm{\alpha}}(\bm{x},0),
\label{eq_simple_expansion_of_time_evolution}
\end{align}
where we take $M$ large enough. Although this simple expression means the repeated application of $(\bm{1} + (T/M)\mathcal{L}^\dagger)$, we will see other types of expressions. 

The basis expansion of \eqref{eq_basis_expansion} and the explicit matrix representation in \eqref{eq_explicit_expression_for_operetor} yield an algorithm based on combinatorics. An explicit example helps us understand the algorithm. Hence, we introduce an illustrative example next.

\subsection{Illustrative example}
\label{subsec_illustrative_example}

\begin{table}
\caption{Terms in the adjoint time-evolution operator for the noisy van der Pol system.}
\label{table_event_van_der_Pol}
\begin{center}
\begin{tabular}{cl}
\hline\hline
No. & Term \\
\hline 
{1} & $(1/2) \nu_{11}^2 \partial_{x_1}^{2}$  \\
{2} & $x_{\mathrm{ini},2} \partial_{x_{1}}$ \\
{3} & $(x_{2}-x_{\mathrm{ini},2}) \partial_{x_{1}}$ \\
{4} & $(1/2) \nu_{22}^{2} \partial_{x_{2}}^{2}$ \\
{5} & $- x_{\mathrm{ini},1} \partial_{x_{2}} - \epsilon \left(x_{\mathrm{ini},1}\right)^{2} x_{\mathrm{ini},2} \partial_{x_{2}} + \epsilon x_{\mathrm{ini},2} \partial_{x_{2}} \quad$ \\
{6} & $- \epsilon \left(x_{\mathrm{ini},1}\right)^{2} (x_{2}-x_{\mathrm{ini},2})  \partial_{x_{2}} + \epsilon (x_{2}-x_{\mathrm{ini},2})  \partial_{x_{2}} $ \\
{7} & $- 2 \epsilon x_{\mathrm{ini},1} x_{\mathrm{ini},2} (x_{1}-x_{\mathrm{ini},1}) \partial_{x_{2}} - (x_{1}-x_{\mathrm{ini},1}) \partial_{x_{2}}$ \\
{8} & $- 2 \epsilon x_{\mathrm{ini},1} (x_{1}-x_{\mathrm{ini},1}) (x_{2}-x_{\mathrm{ini},2})  \partial_{x_{2}}$ \\
{9} & $ - \epsilon x_{\mathrm{ini},2} (x_{1}-x_{\mathrm{ini},1})^{2} \partial_{x_{2}} $ \\
{10} & $- \epsilon (x_{1}-x_{\mathrm{ini},1})^{2} (x_{2}-x_{\mathrm{ini},2})  \partial_{x_{2}} $  \\
\hline\hline
\end{tabular}
\end{center}
\end{table}

\begin{table}
\caption{Coefficients and state-change vectors for Table~\ref{table_event_van_der_Pol}.}
\label{table_event_van_der_Pol2}
\begin{center}
\begin{tabular}{cll}
\hline\hline
No.  & $\displaystyle \gamma_r(\bm{n})$ & $\displaystyle \bm{v}_r$\\
\hline
{1} 
& $ (1/2) \nu_{11}^2 n_1(n_1-1)$
& $[-2,0]$ \\
{2} 
& $x_{\mathrm{ini},2} n_1$
& $[-1,0]$ \\
{3} 
& $n_1$
& $[-1,1]$ \\
{4} 
& $(1/2) \nu_{22}^{2} n_2(n_2-1)$ 
& $[0,-2]$ \\
{5} 
& $- x_{\mathrm{ini},1} n_2 - \epsilon \left(x_{\mathrm{ini},1}\right)^{2} x_{\mathrm{ini},2} n_2 + \epsilon x_{\mathrm{ini},2} n_2 \quad$
& $[0,-1]$ \\
{6} 
& $- \epsilon \left(x_{\mathrm{ini},1}\right)^{2} n_2 + \epsilon n_2$
& $[0,0]$\\
{7} 
& $- 2 \epsilon x_{\mathrm{ini},1} x_{\mathrm{ini},2} n_2 - n_2$
& $[1,-1]$\\
{8} 
& $- 2 \epsilon x_{\mathrm{ini},1} n_2$
& $[1,0]$\\
{9} 
& $ - \epsilon x_{\mathrm{ini},2} n_2$
& $[2,-1]$\\
{10}
& $- \epsilon n_2$
& $[2,0]$\\
\hline\hline
\end{tabular}
\end{center}
\end{table}

The illustrative example is the noisy van der Pol system, which has two stochastic variables:
\begin{align}
\left\{ \begin{array}{l}
dX_1 = X_2 dt + \nu_{11} dW_1(t), \\
dX_2 = \left( \epsilon X_2 \left( 1 - X_1^2 \right) - X_1\right) dt + \nu_{22} dW_2(t),
\end{array}\right.
\label{eq_van_der_pol_equation}
\end{align}
where $\epsilon > 0$, $\nu_{11} > 0$ and $\nu_{22} > 0$. The noisy version of the van der Pol system \cite{van_der_Pol1926} has already been used in filtering \cite{Lakshmivarahan2009} and Koopman operators \cite{Crnjaric-Zic2020}. The adjoint time-evolution operator $\mathcal{L}^\dagger$ is
\begin{align}
\mathcal{L}^\dagger =& x_2 \partial_{x_1} + \left( \epsilon x_2 \left( 1 - x_1^2 \right) - x_1\right) \partial_{x_2}
+ \frac{1}{2} \nu_{11}^2 \partial_{x_1}^2 + \frac{1}{2} \nu_{22}^2 \partial_{x_2}^2.
\end{align}
Introducing $\bm{x}_\mathrm{ini} = [x_{\mathrm{ini},1}, x_{\mathrm{ini},2}]^\mathrm{T}$, we finally obtain the ten events summarized in Table~\ref{table_event_van_der_Pol}, which leads to coefficients and state-change vectors in Table~\ref{table_event_van_der_Pol2}.

\begin{figure}
\begin{center}
\includegraphics[width=48mm]{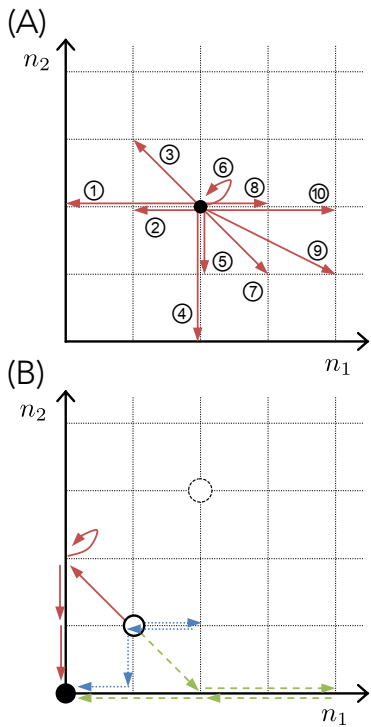}
\end{center}
\caption{\label{fig_combinatorics} (Color online) (A) Ten events for the noisy van der Pol system in Table~\ref{table_event_van_der_Pol2} on a two-dimensional lattice. Each event causes a different process on the lattice. (B) Some examples of walks with $M = 4$ to evaluate the target statistics. All walks finish at $\bm{n} = \bm{0}$ (the filled circle). The initial coordinate corresponds to the target statistics; the solid-empty circle gives $\mathbb{E}[(X_1(T) - x_{\mathrm{ini},1}) (X_2(T) - x_{\mathrm{ini},2}) |\bm{X}(0) = \bm{x}_{\mathrm{ini}}]$, and the dashed-empty circle yields $\mathbb{E}[(X_1(T)-x_{\mathrm{ini},1})^2 (X_2(T) - x_{\mathrm{ini},2})^3|\bm{X}(0) = \bm{x}_{\mathrm{ini}}]$.}
\end{figure}

Figure~\ref{fig_combinatorics}(A) shows ten events for the noisy van der Pol system in Table~\ref{table_event_van_der_Pol}. The coordinate on the two-dimensional lattice corresponds to $\bm{n}$. An action of $\bm{1} + (T/M)\mathcal{L}^\dagger$ in \eqref{eq_simple_expansion_of_time_evolution} gives a `process' corresponding to an event on the lattice. Note that there is a process that does not change the coordinate. We need to enumerate all walks from $\bm{n} = \bm{\alpha}$ to $\bm{n} = \bm{0}$ for the calculation of $P_{\bm{\alpha}}(\bm{0},T)$ in \eqref{eq_zero_coefficient}. Each event causes the weight factor $\gamma_r(\bm{n})$ in \eqref{eq_explicit_expression_for_operetor} that depends on the event number and the coordinate on the lattice. The weights for the walks are used to calculate the target statistics. Figure~\ref{fig_combinatorics}(B) gives some examples of walks with $M = 4$ to evaluate the target statistics. Note that all walks finish at $\bm{n} = \bm{0}$ (the filled circle). The initial coordinate corresponds to the target statistics, i.e., the $\bm{\alpha}$-th moment with the shift of origin. For example, the solid-empty circle gives $\mathbb{E}[(X_1(T) - x_{\mathrm{ini},1}) (X_2(T) - x_{\mathrm{ini},2}) |\bm{X}(0) = \bm{x}_{\mathrm{ini}}]$, and the dashed-empty circle yields $\mathbb{E}[(X_1(T)-x_{\mathrm{ini},1})^2 (X_2(T) - x_{\mathrm{ini},2})^3|\bm{X}(0) = \bm{x}_{\mathrm{ini}}]$.

\subsection{Short summary}
\label{subsec_short_summary}

The algorithm based on combinatorics employs the time-discretization in \eqref{eq_simple_expansion_of_time_evolution}. However, we do not need approximations such as the space discretization and the truncation of equations. As we will see later, a finite $M$ for the time-discretization gives a good estimation for the target statistics. Although the number of the coefficients of the basis expansion in \eqref{eq_basis_expansion} is infinite, only a finite number of $P_{\bm{\alpha}}(\bm{n},T)$ is non-zero for the finite $M$. Hence, there is no need to employ the truncation for $\bm{n}$.

It is time-consuming to enumerate all walks. Hence, an algorithm based on dynamic programming should be employed \cite{Mehlhorn_book}. The basic concept is simple; we calculate the weights of the related coordinates on the lattice instead of the enumeration of walks. The initial coordinate is $\bm{\alpha}$, which stems from the target statistics, i.e., the $\bm{\alpha}$-th moment. An action of $1+(T/M)\mathcal{L}^\dagger$ changes the weights of the related coordinates; after repeating the $M$ processes, the weight on $\bm{n} = 0$ corresponds to the value of the target statistics when $M$ is large enough.

Some previous studies showed that the algorithm based on combinatorics achieves rapid computation compared with the Monte Carlo samplings. For example, the algorithm yielded a speedup of more than $40$ times for the case of the noisy van der Pol system in Ref.~\cite{Ohkubo2022}. Of course, the speedup rate depends on systems and parameters. The most crucial point is the length $M$ of the walk; if $M$ is large, the computational time increases rapidly. Hence, we want to seek the improvement of the algorithm so that a small $M$ is enough; the proposals in the following sections contribute to this aim.

\section{Proposal 1: Approximations of Resolvent}

\subsection{Preparation 1: Simple explicit methods}
\label{subsec_simple_explicit_methods}

In \eqref{eq_simple_expansion_of_time_evolution}, we employ a simple time-splitting for $\exp(\mathcal{L}^\dagger T)$. Of course, a higher-order expansion of the exponential gives a higher-order approximation. We summarize two simple algorithms here:

\begin{itemize}
\item Explicit 1st-order method:
This corresponds to \eqref{eq_simple_expansion_of_time_evolution}:
\begin{align}
\exp \left( \frac{T}{M} \mathcal{L}^\dagger \right) \simeq \bm{1} + \frac{T}{M} \mathcal{L}^\dagger.
\label{eq_explicit_Euler}
\end{align}
The convergence rate is 1st-order for $T/M$.
\item Explicit 2nd-order method: A higher-order expansion for the exponential function is available; with the aid of an analogy with the Heun method, we have
\begin{align}
\exp \left( \frac{T}{M} \mathcal{L}^\dagger \right) \simeq \bm{1} + \frac{T}{M} \mathcal{L}^\dagger + \frac{1}{2} \left(\frac{T}{M} \mathcal{L}^\dagger \right)^2,
\label{eq_Heun}
\end{align}
which yields the 2nd-order convergence rate for $T/M$.
\end{itemize}

Note that the explicit 1st-order method needs `one process' in the time-evolution for $T/M$; only one event occurs in the application of \eqref{eq_explicit_Euler}. By contrast, the explicit 2nd-order method takes at most `two processes' for $T/M$; see the third term in the r.h.s. of \eqref{eq_Heun}. Hence, the computational cost in a time step $T/M$ for the explicit 2nd-order method is larger than that for the explicit 1st-order method.

\subsection{Preparation 2: Resolvent}
\label{subsec_resolvent}

In Ref.~\cite{Ohkubo2021}, a method based on the resolvent of operator $\mathcal{L}^\dagger$ was employed. The method is based on 1st-order approximation of the resolvent. From an analogy with implicit Euler method, we call it the implicit 1st-order method. The following expression of the exponential of the time-evolution operator is employed:
\begin{align}
e^{\mathcal{L}^\dagger T}
= \lim_{M \to \infty} \left[ \left(1 - \frac{T}{M} \mathcal{L}^\dagger \right)^{-1} \right]^M,
\label{eq_resolvent}
\end{align}
where $\left(1 - \frac{T}{M} \mathcal{L}^\dagger \right)^{-1}$ is a resolvent of $\mathcal{L}^\dagger$, apart from a constant factor \cite{Kato_book}.

In Ref.~\cite{Ohkubo2021} and Ref.~\cite{Ohkubo2022}, the following approximation of the matrix elements for the resolvent was introduced:
\begin{align}
&\left[ 1 - \frac{T}{M} \mathcal{L}^\dagger \right]_{\bm{n}' \bm{n}}^{-1}
\simeq
\begin{cases}
\displaystyle 
\frac{1}{1-\frac{T}{M} \left[ \mathcal{L}^\dagger \right]_{\bm{n}\bm{n}}}
\quad \textrm{for } \bm{n}' = \bm{n},
\\
\\
\displaystyle
\frac
{\frac{T}{M} \left[ \mathcal{L}^\dagger \right]_{\bm{n}' \bm{n}}}
{\left( 1-\frac{T}{M} \left[ \mathcal{L}^\dagger \right]_{\bm{n}'\bm{n}'} \right) 
\left( 1-\frac{T}{M} \left[ \mathcal{L}^\dagger \right]_{\bm{n}\bm{n}} \right)} \\
\hspace{42mm} \textrm{otherwise.}
\end{cases}
\label{eq_resolvent_explicit}
\end{align}

\subsection{Usage of second-order approximation}
\label{subsec_CN}

In \eqref{eq_resolvent}, the resolvent was naturally introduced from the definition of the exponential of the operator. There is another approach; formally, we obtain the following expression via the Taylor expansion:
\begin{align}
\exp \left( \frac{T}{M} \mathcal{L}^\dagger \right) = \left[\exp \left( - \frac{T}{M} \mathcal{L}^\dagger \right) \right]^{-1} \simeq \left( \bm{1} - \frac{T}{M} \mathcal{L}^\dagger \right)^{-1}.
\label{eq_implicit_Euler_2}
\end{align}
Note that $\mathcal{L}^\dagger$ generates only a semigroup, and we should interpret the deformation as a formal one.
Hence, from an analogy with the Crank-Nicolson method, we have
\begin{align}
\exp \left( \frac{T}{M} \mathcal{L}^\dagger \right) &= 
\exp \left( \frac{T}{2M} \mathcal{L}^\dagger \right) \exp \left( \frac{T}{2M} \mathcal{L}^\dagger \right)
\nonumber \\
&=
\left[ \exp \left( - \frac{T}{2M} \mathcal{L}^\dagger \right) \right]^{-1} \exp \left( \frac{T}{2M} \mathcal{L}^\dagger \right)
\nonumber \\
&\simeq 
\left( \bm{1} - \frac{T}{2M} \mathcal{L}^\dagger \right)^{-1}
\left( \bm{1} + \frac{T}{2M} \mathcal{L}^\dagger \right).
\label{eq_proposed_method}
\end{align}
We here call it the implicit 2nd-order method. The present work focuses on algorithms based on combinatorics, which need explicit expressions for the resolvent. However, we cannot use the approximations in \eqref{eq_resolvent_explicit} because they do not include the 2nd-order contribution adequately. Next, the 2nd-order approximation is explicitly given.

\subsection{Second-order approximation of resolvent}
\label{subsec_inverse_CN}

The aim here is to obtain an explicit matrix expression for the resolvent:
\begin{align}
\mathcal{C} = \left(I - \frac{T}{2M} \mathcal{L}^\dagger \right)^{-1},
\label{eq_CN_inverse}
\end{align}
which is necessary to apply the implicit 2nd-order method for the time-evolution of the backward Kolmogorov equation. 

 We expect \eqref{eq_proposed_method} has the 2nd-order convergence in time. This fact suggests that only the 2nd-order approximation is enough for the inverse part. The appendix gives a detailed discussion of the derivation; we here see only the final expressions. Employing the basis expansion in \eqref{eq_basis_expansion}, the matrix representation of $\mathcal{C}$ is approximated as follows:
\begin{align}
\left[\mathcal{C}\right]_{\bm{n}\bm{n}} \simeq \left(1 - h \left[\mathcal{L}^\dagger\right]_{\bm{n}\bm{n}}
 - h^2 \sum_{\widetilde{\bm{n}} \neq \bm{n}} \left[\mathcal{L}^\dagger\right]_{\bm{n} \widetilde{\bm{n}}} \left[\mathcal{L}^\dagger\right]_{\widetilde{\bm{n}}\bm{n}}  \right)^{-1}
\label{eq_CN_inverse_1}
\end{align}
and
\begin{align}
\left[\mathcal{C}\right]_{\bm{n}'\bm{n}} 
\simeq&
\frac{h \left[\mathcal{L}^\dagger\right]_{\bm{n}'\bm{n}}}{
\left(\displaystyle 1 -  h \left[\mathcal{L}^\dagger\right]_{\bm{n}'\bm{n}'} 
 \right)
\left(\displaystyle 1 -  h \left[\mathcal{L}^\dagger\right]_{\bm{n}\bm{n}} 
 \right)
} \nonumber \\
& + \sum_{\widetilde{\bm{n}} \neq \bm{n}', \bm{n}}
h^2 \left[\mathcal{L}^\dagger\right]_{\bm{n}'\widetilde{\bm{n}}} \left[\mathcal{L}^\dagger\right]_{\widetilde{\bm{n}}\bm{n}},
\label{eq_CN_inverse_2}
\end{align}
where $h = T / (2M)$. 

Equations in \eqref{eq_CN_inverse_1} and \eqref{eq_CN_inverse_2} give approximations up to the 2nd-order in $h$, which is enough to obtain the 2nd-order convergence for the implicit 2nd-order method. Note that \eqref{eq_CN_inverse_1} and \eqref{eq_CN_inverse_2} contain only `local' movements that stem from two successive processes. Furthermore, the backward Kolmogorov equation leads to a kind of sparse structure for $\mathcal{L}^\dagger$ for the expansion in \eqref{eq_basis_expansion}, as shown in Table~\ref{table_event_van_der_Pol2}. Hence, the computational cost is not high, and algorithms based on combinatorics are available. 

There is a comment on the number of processes. \eqref{eq_proposed_method} needs at most `three processes' for $T/M$; one might consider that the computational cost is high. However, as we will see later, the proposed method gives a sufficiently accurate value with a small $M$.

In the following, we use two examples in order to check the convergence property of the proposed method.

\subsection{Numerical experiments 1: Ornstein-Uhlenbeck process}
\label{subsec_experiments_OU}

The first example is the famous Ornstein-Uhlenbeck process  \cite{Gardiner_book}:
\begin{align}
dX = - \gamma X dt + \sigma dW(t),
\end{align}
where $\gamma > 0$ and $\sigma > 0$. The adjoint time-evolution operator $\mathcal{L}^\dagger$ for the backward Kolmogorov equation is
\begin{align}
\mathcal{L}^\dagger = - \gamma x \frac{\partial}{\partial x} 
+ \frac{\sigma^2}{2} \frac{\partial^2}{\partial x^2}.
\end{align}

\begin{table}[t]
\caption{Coefficients and state-change vectors for the Ornstein-Uhlenbeck process.}
\label{table_event_OU}
\begin{center}
\begin{tabular}{lll}
\hline\hline
Term & $\displaystyle \gamma_r(\bm{n})$ & $\displaystyle \bm{v}_r$\\
\hline
$- \gamma (x-x_\mathrm{ini}) \partial_{x}$ \quad
& $ - \gamma n$ \quad
& $[0]$ \\
$- \gamma x_\mathrm{ini} \partial_{x}$  \quad
& $- \gamma n$ \quad
& $[-1]$ \\
$(\sigma^2 / 2) \partial_{x}^2$ \quad
& $( \sigma^2 / 2) n(n-1)$ \quad
& $[-2]$ \quad \\
\hline\hline
\end{tabular}
\end{center}
\end{table}

In order to employ the basis expansion in \eqref{eq_basis_expansion}, we rewrite $\mathcal{L}^\dagger$ as follows:
\begin{align}
\mathcal{L}^\dagger &= - \gamma (x-x_\mathrm{ini}) \frac{\partial}{\partial x} 
- \gamma x_\mathrm{ini} \frac{\partial}{\partial x} 
+ \frac{\sigma^2}{2} \frac{\partial^2}{\partial x^2}.
\label{eq_OU_adjoint_operator}
\end{align}
The three terms in \eqref{eq_OU_adjoint_operator} correspond to three events shown in Table~\ref{table_event_OU}.

The following expectation values are obtained using the established analytical solutions of the Ornstein-Uhlenbeck process \cite{Gardiner_book}:
\begin{align}
&\mathbb{E}\left[ X_T - x_\mathrm{ini} | X_0 = x_\mathrm{ini}\right] = x_\mathrm{ini} e^{-\gamma T} - x_\mathrm{ini},
\label{eq_OU_1st_moment} \\
&\mathbb{E}\left[ (X_T - x_\mathrm{ini})^2 | X_0 = x_\mathrm{ini} \right] \nonumber \\
&\, = \frac{\sigma^2}{2\gamma} \left( 1- e^{-2\gamma T}\right)
+ x_\mathrm{ini}^2 e^{-2\gamma T} 
- 2 x_\mathrm{ini}^2 e^{-\gamma T}
+ x_\mathrm{ini}^2,
\label{eq_OU_2nd_moment}
\end{align}
where the initial coordinate is $x_\mathrm{ini}$. In the numerical experiments below, we use $\gamma = 1.0$, $\sigma = 0.5$ and $x_\mathrm{ini} = 1.0$.

\begin{figure}
\begin{center}
\includegraphics[width=90mm]{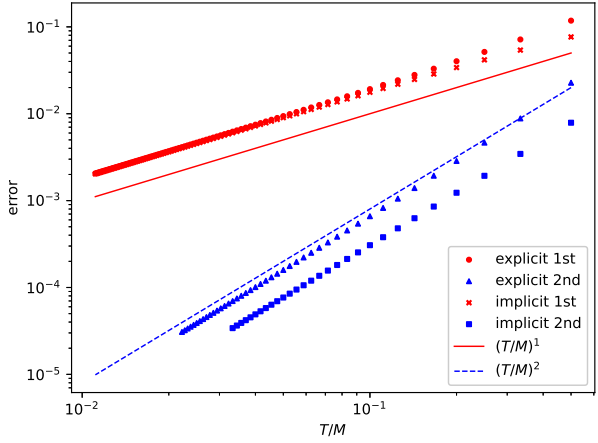}
\end{center}
\caption{\label{fig_OU_1st} (Color online) Errors from the exact solutions for the 1st-order moment in the Ornstein-Uhlenbeck process for $T = 1.0$. Circle, triangle, cross, and square markers correspond to the explicit 1st-order method, the explicit 2nd-order method, the implicit 1st-order method, and the implicit 2nd-order method, respectively.}
\end{figure}

Here, we check the convergence property for $T/M$. For the Ornstein-Uhlenbeck case, the evaluation of the error from the exact solutions is easy because there are the analytical solutions. For example, the 1st-order moment is evaluated for several different methods for $T = 1.0$. Figure~\ref{fig_OU_1st} shows the results in which circle, triangle, cross, and square markers correspond to the explicit 1st-order method, the explicit 2nd-order method, the implicit 1st-order method, and the implicit 2nd-order method, respectively. The results in Figure~\ref{fig_OU_1st} indicate that The expected convergence rates are obtained. Several other experiments show the same convergence rate for other cases.

\subsection{Numerical experiments 2: noisy van der Pol process}
\label{subsec_experiments_van_der_Pol}

The second example is the noisy van der Pol process in \eqref{eq_van_der_pol_equation}. In the numerical experiments below, we set $\epsilon = 1.0$, $\nu_{11} = 0.5$, and $\nu_{22} = 0.5$. The initial coordinate is set as $x_{\mathrm{ini},1} = 0.5$ and $x_{\mathrm{ini},2} = 1.0$.

There is a comment on the approximation for the resolvent. In the proposed method, we employ the 2nd-order approximation for the resolvent (inverse matrix) in \eqref{eq_CN_inverse_1} and \eqref{eq_CN_inverse_2}. If the 1st-order approximation for the resolvent in \eqref{eq_resolvent_explicit} is employed for the Crank-Nicolson method, the convergence behavior shows the 1st-order one. Hence, the proposed approximation in \eqref{eq_CN_inverse_1} and \eqref{eq_CN_inverse_2} is crucial to obtain the convergence behavior of the proposed method.

\begin{figure}
\begin{center}
\includegraphics[width=90mm]{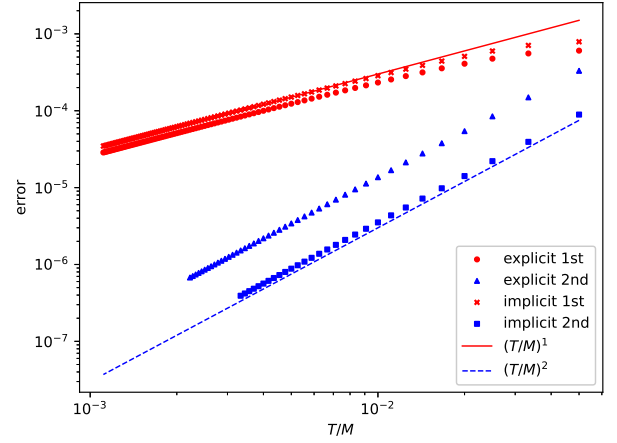}%
\end{center}
\caption{\label{fig_van_der_Pol_cross} (Color online) Errors from the numerically obtained exact solutions for $\mathbb{E}[(X_1(T)-x_{\mathrm{ini},1}) (X_2(T) -x_{\mathrm{ini},2})| \bm{X}(0) = \bm{x}_\mathrm{ini}]$ in the noisy van der Pol system for $T = 0.1$. 
}
\end{figure}

Different from the Ornstein-Uhlenbeck case, we cannot obtain exact analytical solutions to the moments for the noisy van der Pol system. Hence, the 4th-order Runge-Kutta method with $\Delta t = 10^{-10}$ is employed for $T = 0.1$, and the coupled ordinary differential equations for $P_{\bm{\alpha}}(\bm{n},t)$ are numerically solved. The finite cutoff of the state should be employed, and $n_1 , n_2 < 15$ is used here; note that a larger cutoff does not change the final numerical values. It is possible to see the calculated values from the coupled ordinary differential equations as numerically obtained exact solutions because of the short time interval $\Delta t = 10^{-10}$. 

Figure~\ref{fig_van_der_Pol_cross} shows the error behavior for $\mathbb{E}[(X_1(T)-x_{\mathrm{ini},1}) (X_2(T) -x_{\mathrm{ini},2})| \bm{X}(0) = \bm{x}_\mathrm{ini}]$, and the obtained convergence properties are the same as the expected one.

\subsection{Comments for computational time}
\label{subsec_result_computationl_time}

The main aim of the present paper is the proposition of the approximations of resolvent. Hence, the code is not optimized. However, the computational time is briefly mentioned here. 

For the noisy van der Pol system, the coupled ordinary differential equations for $P_{\bm{\alpha}}(\bm{n},t)$ give the numerical evaluation of moments. However, the number of equations is infinite in general, and we need some finite cutoff, as stated above. For the case with $T = 0.1$ and the parameters above, the cutoff with $n_1, n_2 < 15$ was used. Despite the small cutoff size, $\Delta t = 10^{-10}$ needs high computational time. The computational time was about $40$ hours with Intel(R) Xeon(R) Gold 6128 CPU, 3.40GHz. 

The experiments for the proposed algorithms were performed on MacBookAir with an M1 processor. Note that there is no need to employ an explicit cutoff for $\bm{n}$ because the algorithm based on combinatorics considers all walks. Of course, the length of the walks is limited by the parameter $M$. In the implicit 2nd-order method, a case with $M = 30$ takes about $24$ seconds; that with $M=20$ needs about $3$ seconds, and that with $M=10$ is finished in less than $0.2$ seconds.

The increase of $M$ needs exponential increases in computational time. Hence, it is preferable to use small $M$. However, as shown in Figs.~\ref{fig_OU_1st} and \ref{fig_van_der_Pol_cross}, small $M$ cases are not enough to obtain accurate results; the convergence behaviors are not fast. Next, we try to improve the convergence with the aid of extrapolation.

\section{Proposal 2: Extrapolation}
\label{sec_extrapolation}

We usually perform the algorithms based on combinatorics with a parameter $M$. Only one estimation is not enough to judge whether the estimated value is converged. Hence, we repeat the algorithms a few times for different values of $M$.

Note that we checked the convergent properties of the algorithms in Figs.~\ref{fig_OU_1st} and \ref{fig_van_der_Pol_cross}. These results suggest that the extrapolation with $T/M \to 0$ is available to obtain more accurate estimations.

\subsection{First-order cases}

Consider two calculations for two different $M$ cases, $M^{(i)}$ ($i = 1, 2$). Let $\widehat{m}^{(i)}$ be the evaluated value of the target statistic for the case with $M^{(i)}$. When the true value is $\widetilde{m}$, the 1st-order convergence suggests the following relations;
\begin{align}
&\left| \widehat{m}^{(1)} - \widetilde{m} \right| = C \left( \frac{T}{M^{(1)}}\right),  \qquad
\left| \widehat{m}^{(2)} - \widetilde{m} \right| = C \left( \frac{T}{M^{(2)}}\right),
\end{align}
where $C$ is a constant.

If $\widehat{m}^{(1)} > \widetilde{m}$ and $\widehat{m}^{(2)} > \widetilde{m}$, the following estimation is derived:
\begin{align}
\widetilde{m} = \widehat{m}^{(2)} - \frac{M^{(1)} \left( \widehat{m}^{(1)} - \widehat{m}^{(2)}\right) }{(M^{(2)} - M^{(1)})}.
\label{eq_first_order_extrapolation}
\end{align}
For $\widehat{m}^{(1)} < \widetilde{m}$ and $\widehat{m}^{(2)} < \widetilde{m}$, the same estimation is derived. Preliminary numerical experiments show that the behavior of $\widehat{m}^{(i)}$ is similar for small $T/M^{(i)}$. Hence, \eqref{eq_first_order_extrapolation} is available.

In practice, we perform the algorithm with $M^{(1)} = M-1$ and obtain the final result $m^{(1)}$. Then, we employ $M^{(2)} = M$, and the algorithm gives $m^{(2)}$. Finally, using \eqref{eq_first_order_extrapolation}, we estimate $\widetilde{m}$ as the new expected value for the target statistic.

\subsection{Second-order cases}

For the case with the 2nd-order convergence, we have
\begin{align}
&\left| \widehat{m}^{(1)} - \widetilde{m} \right| = C \left( \frac{T}{M^{(1)}}\right)^2,  \qquad
\left| \widehat{m}^{(2)} - \widetilde{m} \right| = C \left( \frac{T}{M^{(2)}}\right)^2.
\end{align}
Hence, the following extrapolation method is available as similar to the 1st-order case:
\begin{align}
\widetilde{m} = \widehat{m}^{(2)} -  \frac{\left( M^{(1)} \right)^2 \left( \widehat{m}^{(1)} - \widehat{m}^{(2)}\right)}{ \left( \left(M^{(2)}\right)^2 - \left(M^{(1)}\right)^2 \right)}. 
\label{eq_second_order_extrapolation}
\end{align}

\subsection{Numerical checks}

\begin{figure}
\includegraphics[width=90mm]{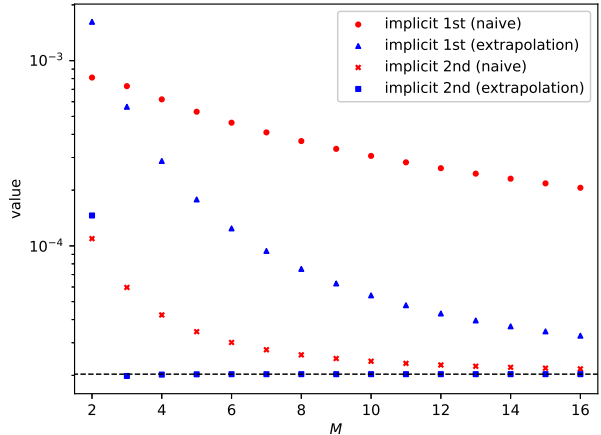}%
\caption{\label{fig_van_der_Pol_extrapolation} (Color online) Extrapolation results for the noisy van der Pol system. The evaluated statistics is $\mathbb{E}[(X_1(T)-x_{\mathrm{ini},1}) (X_2(T) -x_{\mathrm{ini},2})| \bm{X}(0) = \bm{x}_\mathrm{ini}]$, and $T = 0.1$. Circle, triangle, cross, and square markers correspond to the raw values obtained from the implicit 1st-order method, the extrapolation in \eqref{eq_first_order_extrapolation}, the raw values obtained from the implicit 2nd-order method, and the extrapolation in \eqref{eq_second_order_extrapolation}, respectively. The dashed line is the numerically obtained exact solution, $2.030 \times 10^{-5}$.
}
\end{figure}

Figure~\ref{fig_van_der_Pol_extrapolation} shows the estimated values for $\mathbb{E}[(X_1(T)-x_{\mathrm{ini},1}) (X_2(T) -x_{\mathrm{ini},2})| \bm{X}(0) = \bm{x}_\mathrm{ini}]$ for the same settings in Sec.~\ref{subsec_experiments_van_der_Pol}; circle, triangle, cross, and square markers correspond to the raw values obtained from the implicit 1st-order method, the extrapolation in \eqref{eq_first_order_extrapolation}, the raw values obtained from the implicit 2nd-order method, and the extrapolation in \eqref{eq_second_order_extrapolation}, respectively. The exact solution is estimated numerically by the 4th-order Runge-Kutta method with $\Delta t = 10^{-10}$; the value is $2.030 \times 10^{-5}$. We can see that the extrapolations give rapid convergence compared with the raw values. Since a small $M$ is enough to obtain accurate results, the extrapolations reduce the computational cost.

\subsection{Higher-order moment}

\begin{figure}
\includegraphics[width=90mm]{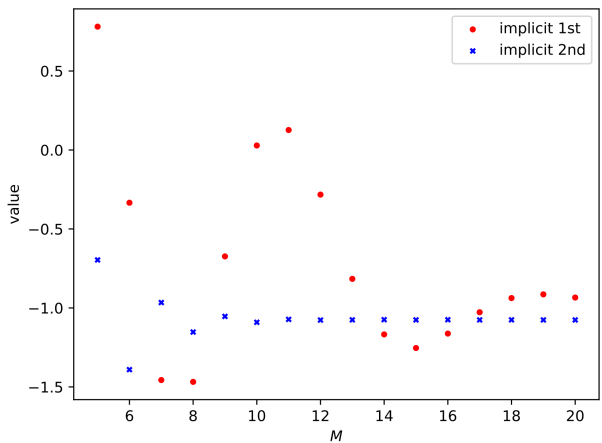}%
\caption{\label{fig_van_der_Pol_higher_order} (Color online) An example for higher-order moment cases. The evaluated statistics is $\mathbb{E}[(X_1(T)-x_{\mathrm{ini},1})^3 (X_2(T) - x_{\mathrm{ini},2})^3| \bm{X}(0) = \bm{x}_\mathrm{ini}]$, where $T=1.0$ and $\bm{x}_\mathrm{ini} = (0.5,1.0)$. Circle and cross markers correspond to the values obtained from the implicit 1st-order method and the implicit 2nd-order method, respectively.
}
\end{figure}

Finally, we see an example for higher-order moment cases, which will show one of the benefits of the higher-order approximation. Figure~\ref{fig_van_der_Pol_higher_order} shows the numerical results for $\mathbb{E}[(X_1(T)-x_{\mathrm{ini},1})^3 (X_2(T) - x_{\mathrm{ini},2})^3| \bm{X}(0) = \bm{x}_\mathrm{ini}]$, where $T=1.0$ and $\bm{x}_\mathrm{ini} = (0.5,1.0)$; circle and cross markers correspond to the values obtained from the implicit 1st-order method and the implicit 2nd-order method, respectively. Although evaluations for higher-order moments could be unstable generally, as shown in the 1st-order method in Fig.~\ref{fig_van_der_Pol_higher_order}, the implicit 2nd-order one gives more stable behavior in the small $M$ regions. Of course, we need further investigations and modifications in the future, which will include higher-order approximations of resolvents. However, the proposed 2nd-order method will be a hopeful candidate to evaluate higher-order moments that could be important to seek large deviations or nonequilibrium characteristics of systems.

\section{Conclusions}

The aim of the present paper is not to obtain the solution of the backward Kolmogorov equation on the whole domain. Hence, algorithms based on combinatorics are available. To construct the proposed algorithm with 2nd-order convergence, the explicit matrix representation for the resolvent up to the 2nd-order is derived; the matrix elements are calculated only from local processes. Furthermore, the convergence property naturally gives the extrapolation method, which reduces the computational costs.

Although the power of the local approximation of the resolvent is demonstrated in the present paper, there are some remaining studies. For example, we should discuss the stability of numerical algorithms in the future; the approximation of the resolvent (inverse matrix) in \eqref{eq_CN_inverse_1} and \eqref{eq_CN_inverse_2} would vary the stability property. In addition, although the proposed algorithm based on combinatorics works well when the number of variables is small, the computation time grows exponentially when the number of variables is large. Hence, we need further approximations. Note that a similar discussion based on combinatorics was applied to a partial differential equation in control theory \cite{Ohkubo2021b}. While the aim was not to evaluate statistics, we can employ the discussions in the present paper to such different research fields. Hence, approximation studies for algorithms based on combinatorics are also important for future work.

\appendix

\begin{widetext}

\section{Approximation of the resolvent}
\label{sec_appendix_inverse}

We need to evaluate the matrix representation of the resolvent, $\left( \bm{1} - (T/2M) \mathcal{L}^\dagger \right)^{-1}$ for the implicit 2nd-order method. However, the original $\mathcal{L}^\dagger$ is infinite-dimensional for the basis expansion of \eqref{eq_basis_expansion}. Hence, the operator $\mathcal{L}^\dagger$ gives an infinite matrix representation, and it is not straightforward to obtain the inverse matrix. We here explain how to derive the expressions in \eqref{eq_CN_inverse_1} and \eqref{eq_CN_inverse_2} with the aid of finite cases.

Using the explicit basis expansion of \eqref{eq_basis_expansion}, we obtain a corresponding matrix representation; we denote the matrix element as $[\ell_{ij}]$. Let $D^{(n)}$ be an $n \times n $ matrix constructed by $\left( \bm{1} - (T/2M) \mathcal{L}^\dagger \right)^{-1}$:
\begin{align}
D^{(n)} &= \begin{bmatrix}
d_{11} & d_{12} & d_{13} & \cdots & d_{1n}\\
d_{21} & d_{22} & d_{23} & \cdots & d_{2n} \\
d_{31} & d_{32} & d_{33} & \cdots & d_{3n} \\
\vdots & \vdots & \vdots & \ddots & \vdots \\
d_{n1} & d_{n2} & d_{n3} & \cdots & d_{nn} \\
\end{bmatrix} 
=
\begin{bmatrix}
1 - h \ell_{11} & - h \ell_{12} & - h \ell_{13} & \cdots & - h \ell_{1n} \\
- h \ell_{21} & 1 - h \ell_{22} & - h \ell_{23} & \cdots & - h \ell_{2n} \\
- h \ell_{31} & - h \ell_{32} & 1 - h \ell_{33} & \cdots & - h \ell_{3n} \\
\vdots & \vdots & \vdots & \ddots & \vdots \\
- h \ell_{n1} & - h \ell_{n2} & - h \ell_{n3} & \cdots & 1 - h \ell_{nn} \\
\end{bmatrix},
\end{align}
where $h = T/(2M)$. Although we here discuss the finite matrix $D^{(n)}$ for simplicity, it is possible to apply the following discussions formally for $n \to \infty$ case.

Here, we employ the following definition of the determinant of a matrix. Using $\mathcal{I}^{(n)} = \{1,2,\dots,n\}$, let $\mathcal{S}^{(n)}$ be a symmetric group of $\mathcal{I}^{(n)}$. Then, using the signature of the permutation for $\sigma \in \mathcal{S}^{(n)}$, $\mathrm{sgn}(\sigma)$, the determinant is defined as \cite{Beilina_book}
\begin{align}
\mathrm{det}\left(D^{(n)}\right) = \sum_{\sigma \in \mathcal{S}^{(n)}} \mathrm{sgn}(\sigma) \prod_{k=1}^n d_{k \sigma(k)},
\label{eq_appendix_definition_det}
\end{align}
where $\sigma(k)$ is the value in the $k$-th position after the reordering $\sigma$.

Let $C^{(n)}$ be the inverse matrix of $D^{(n)}$; $C^{(n)} = \left(D^{(n)}\right)^{-1}$. The value of $h$ is small, and hence we now discuss the expansion of $C^{(n)}$ with respect to $h$. Using the determinant and adjugate, $C^{(n)}$ is calculated as follows \cite{Beilina_book}:
\begin{align}
C^{(n)} = \frac{1}{\mathrm{det}(D^{(n)})}\mathrm{adj}\left(D^{(n)}\right),
\label{eq_appendix_definition_inverse}
\end{align}
where
\begin{align}
\left[\mathrm{adj}\left(D^{(n)}\right)\right]_{ij}
= (-1)^{i+j}
 \mathrm{det}
\begin{pmatrix}
d_{11} & \cdots & d_{1(i-1)} & d_{1(i+1)} & \cdots & d_{1n} \\
\vdots & \ddots & \vdots & \vdots & \ddots & \vdots \\
d_{(j-1)1} & \cdots & d_{(j-1)(i-1)} & d_{(j-1)(i+1)} & \cdots & d_{(j-1)n} \\
d_{(j+1)1} & \cdots & d_{(j+1)(i-1)} & d_{(j+1)(i+1)} & \cdots & d_{(j+1)n} \\
\vdots & \ddots & \vdots & \vdots & \ddots & \vdots \\
d_{n1} & \cdots & d_{n(i-1)} & d_{n(i+1)} & \cdots & d_{nn} 
\end{pmatrix}.
\label{eq_appendix_adj}
\end{align}
The adjugate is defined with the determinant of the $(n-1)\times(n-1)$ matrix which results from deleting row $j$ and column $i$ of $D^{(n)}$.

First, let us discuss the diagonal elements. Here, we introduce a set $\mathcal{I}_{\setminus i}^{(n)} = \{1,2,\dots,i-1,i+1,\dots,n\}$, and let $\mathcal{S}_{\setminus i}^{(n)}$ be a symmetric group of $\mathcal{I}_{\setminus i}^{(n)}$. Then, employing the definition \eqref{eq_appendix_definition_det} or the cofactor expansion of the determinant, we have
\begin{align}
\left[\mathrm{adj}\left(D^{(n)}\right)\right]_{ii}
= \prod_{\substack{k = 1\\(k \neq i)}}^n d_{kk}
 - \sum_{\sigma \in \mathcal{S}_{\setminus i}^{(n)}} \left( \prod_{l=1}^{n-3} d_{\sigma(l) \sigma(l)} \right) d_{\sigma(n-2)\sigma(n-1)}d_{\sigma(n-1)\sigma(n-2)}
+ \mathcal{O}(h^3).
\label{eq_appendix_adj_diag_result}
\end{align}
Note that all the non-diagonal elements, $d_{ij} (i \neq j)$, are $\mathcal{O}(h)$. Hence, the existence of the non-diagonal elements immediately yields higher-order terms. \eqref{eq_appendix_definition_det} implies that we should consider many combinations of $\{d_{kl}\}$. However, there is no 1st-order term because of the permutation in \eqref{eq_appendix_definition_det}; note that \eqref{eq_appendix_adj} indicates the absence of terms related to the $i$-th element. Only the second term in the r.h.s. in \eqref{eq_appendix_adj_diag_result} is $\mathcal{O}(h^2)$.

For the denominator of the r.h.s. in \eqref{eq_appendix_definition_inverse}, we have
\begin{align}
\mathrm{det}\left(D^{(n)}\right)
= \prod_{k=1}^n d_{kk} 
- \sum_{\sigma \in \mathcal{S}^{(n)}} \left( \prod_{l=1}^{n-2} d_{\sigma(l) \sigma(l)} \right) d_{\sigma(n-1)\sigma(n)}d_{\sigma(n)\sigma(n-1)} 
+ \mathcal{O}(h^3).
\label{eq_appendix_det_result}
\end{align}
Hence, \eqref{eq_appendix_definition_inverse}, \eqref{eq_appendix_adj_diag_result} and \eqref{eq_appendix_det_result} yield the diagonal elements for $C^{(n)}$ as
\begin{align}
\left[C^{(n)}\right]_{ii} 
&= \frac{1 - \sum_{\sigma \in \mathcal{S}_{\setminus i}^{(n)}} d_{\sigma(n-2)\sigma(n-1)}d_{\sigma(n-1)\sigma(n-2)} + \mathcal{O}(h^3)}
{d_{ii} - \sum_{\sigma \in \mathcal{S}^{(n)}} d_{\sigma(n-1)\sigma(n)}d_{\sigma(n)\sigma(n-1)}  + \mathcal{O}(h^3)}.
\label{eq_appendix_cii_pre}
\end{align}
Here, note that $d_{ii} = 1 - h \ell_{ii}$ and the order of second term in the denominator is $h^2$. 
In the derivation, we used
\begin{align*}
\sum_{\sigma \in \mathcal{S}^{(n)}} \frac{d_{\sigma(n-1)\sigma(n)}d_{\sigma(n)\sigma(n-1)}}{d_{\sigma(n-1)\sigma(n-1)}d_{\sigma(n)\sigma(n)}} 
= 
\sum_{\sigma \in \mathcal{S}^{(n)}} \frac{d_{\sigma(n-1)\sigma(n)}d_{\sigma(n)\sigma(n-1)}}{1 + \mathcal{O}(h)} 
= 
\sum_{\sigma \in \mathcal{S}^{(n)}} d_{\sigma(n-1)\sigma(n)}d_{\sigma(n)\sigma(n-1)} +  \mathcal{O}(h^3),
\end{align*}
because $d_{ij} (i \neq j)$ is $\mathcal{O}(h)$.
Hence, the second term in the denominator of \eqref{eq_appendix_cii_pre} is derived. The same discussion is applied to the second term in the numerator of \eqref{eq_appendix_cii_pre}.

Then, the Taylor expansion, $(1-x)^{-1} = \sum_{n=0}^{\infty} x^n$, leads to
\begin{align}
&\frac{1 - \sum_{\sigma \in {\mathcal{S}_{\setminus i}}^{(n)}} d_{\sigma(n-2)\sigma(n-1)}d_{\sigma(n-1)\sigma(n-2)} + \mathcal{O}(h^3)}
{d_{ii} - \sum_{\sigma \in {\mathcal{S}}^{(n)}} d_{\sigma(n-1)\sigma(n)}d_{\sigma(n)\sigma(n-1)}  + \mathcal{O}(h^3)} \nonumber \\
&=
\left[
\left(d_{ii} - \sum_{\sigma \in \mathcal{S}^{(n)}} d_{\sigma(n-1)\sigma(n)}d_{\sigma(n)\sigma(n-1)}  + \mathcal{O}(h^3)\right) 
\left(1 - \sum_{\sigma \in {\mathcal{S}_{\setminus i}}^{(n)}}  d_{\sigma(n-2)\sigma(n-1)}d_{\sigma(n-1)\sigma(n-2)}\right)^{-1}
\right]^{-1} + \mathcal{O}(h^3) 
 \nonumber \\
&=
\left[ \left(d_{ii} - \sum_{\sigma \in \mathcal{S}^{(n)}} d_{\sigma(n-1)\sigma(n)}d_{\sigma(n)\sigma(n-1)}  + \mathcal{O}(h^3)\right) 
 \left(1 + \sum_{\sigma \in {\mathcal{S}_{\setminus i}}^{(n)}}  d_{\sigma(n-2)\sigma(n-1)}d_{\sigma(n-1)\sigma(n-2)} + \mathcal{O}(h^3)\right)
\right]^{-1} + \mathcal{O}(h^3) \nonumber \\
&=
\left[ d_{ii} - \sum_{\sigma \in \mathcal{S}^{(n)}} d_{\sigma(n-1)\sigma(n)}d_{\sigma(n)\sigma(n-1)} 
+\left(1-h \ell_{ii}\right) \sum_{\sigma \in {\mathcal{S}_{\setminus i}}^{(n)}}  d_{\sigma(n-2)\sigma(n-1)}d_{\sigma(n-1)\sigma(n-2)} + \mathcal{O}(h^3)
\right]^{-1}  + \mathcal{O}(h^3) \nonumber \\
&=
\left[ d_{ii} - \sum_{\sigma \in \mathcal{S}^{(n)}} d_{\sigma(n-1)\sigma(n)}d_{\sigma(n)\sigma(n-1)} 
+ \sum_{\sigma \in {\mathcal{S}_{\setminus i}}^{(n)}}  d_{\sigma(n-2)\sigma(n-1)}d_{\sigma(n-1)\sigma(n-2)} + \mathcal{O}(h^3)
\right]^{-1} + \mathcal{O}(h^3) \nonumber \\
&=
\frac{1}
{
d_{ii} - \sum_{k \neq i} d_{ik}d_{ki} + \mathcal{O}(h^3)
} + \mathcal{O}(h^3).
\end{align}
Therefore, we obtain
\begin{align}
\left[C^{(n)}\right]_{ii} 
&= \frac{1}{d_{ii} - \sum_{k \neq i} d_{ik}d_{ki}} + \mathcal{O}(h^3).
\label{eq_appendix_cii}
\end{align}
This corresponds to the matrix expression of \eqref{eq_CN_inverse_1}.

For the non-diagonal cases, the following definition of the determinant is useful \cite{Beilina_book}:
\begin{align}
\mathrm{det}\left(D^{(n)}\right) 
=& \sum_{i=1}^n (-1)^{i+j} d_{ij}
\mathrm{det}
\begin{pmatrix}
d_{11} & \cdots & d_{1(j-1)} & d_{1(j+1)} & \cdots & d_{1n} \\
\vdots & \vdots & \vdots & \vdots & \vdots & \vdots \\
d_{(i-1)1} & \cdots & d_{(i-1)(j-1)} & d_{(i-1)(j+1)} & \cdots & d_{(i-1)n} \\
d_{(i+1)1} & \cdots & d_{(i+1)(j-1)} & d_{(i+1)(j+1)} & \cdots & d_{(i+1)n} \\
\vdots & \vdots & \vdots & \vdots & \vdots & \vdots \\
d_{n1} & \cdots & d_{n(j-1)} & d_{n(j+1)} & \cdots & d_{nn} 
\end{pmatrix}
\label{eq_appendix_det_2} \\
=& \sum_{j=1}^n (-1)^{i+j} d_{ij} 
 \mathrm{det}
\begin{pmatrix}
d_{11} & \cdots & d_{1(j-1)} & d_{1(j+1)} & \cdots & d_{1n} \\
\vdots & \vdots & \vdots & \vdots & \vdots & \vdots \\
d_{(i-1)1} & \cdots & d_{(i-1)(j-1)} & d_{(i-1)(j+1)} & \cdots & d_{(i-1)n} \\
d_{(i+1)1} & \cdots & d_{(i+1)(j-1)} & d_{(i+1)(j+1)} & \cdots & d_{(i+1)n} \\
\vdots & \vdots & \vdots & \vdots & \vdots & \vdots \\
d_{n1} & \cdots & d_{n(j-1)} & d_{n(j+1)} & \cdots & d_{nn} 
\end{pmatrix}.
\label{eq_appendix_det_3}
\end{align}
For the numerator of the r.h.s. in \eqref{eq_appendix_definition_inverse}, we obtain the following expression:
\begin{align}
\left[\mathrm{adj}\left(D^{(n)}\right)\right]_{ij} 
&= 
- d_{ij} \left( \prod_{\substack{l=1\\l\neq i,j}}^n d_{ll} \right)
+ \sum_{\substack{k=1\\k\neq i,j}}^n \left[ d_{ik}d_{kj} \left( \prod_{\substack{l=1\\l\neq i,j,k}}^n d_{ll} \right) \right]
+ \mathcal{O}(h^3).
\label{eq_appendix_adj_non_diag_result}
\end{align}
Note that the definitions of the adjugate in \eqref{eq_appendix_adj} and the determinants in \eqref{eq_appendix_det_2} and \eqref{eq_appendix_det_3} indicate the lacks of terms related the $i$-th and the $j$-th elements, except for $d_{ij}$. Hence, in the determinant, we do not have $d_{ii}$ and $d_{jj}$. Employing the permutation in \eqref{eq_appendix_definition_det}, we verify the first and second terms in \eqref{eq_appendix_adj_non_diag_result}.

As for the determinant in the denominator, we again employ \eqref{eq_appendix_det_result}; since the second term in the r.h.s. of \eqref{eq_appendix_det_result} is $\mathcal{O}(h^2)$, the following expression is available:
\begin{align}
\mathrm{det}\left(D^{(n)}\right)
&= \prod_{k=1}^n d_{kk} + \mathcal{O}(h^2).
\label{eq_appendix_det_result_2}
\end{align}
Hence, from \eqref{eq_appendix_definition_inverse}, \eqref{eq_appendix_adj_non_diag_result}, and \eqref{eq_appendix_det_result_2}, we have
\begin{align}
\left[C^{(n)}\right]_{ij} 
&= - \frac{d_{ij}}{d_{ii} d_{jj}}
+ \sum_{\substack{k=1\\k\neq i,j}}^n \frac{d_{ik}d_{kj}}{d_{ii}d_{kk}d_{jj}}
+ \mathcal{O}(h^3).
\label{eq_appendix_cij_pre}
\end{align}
The computational cost for the denominator in the second term in the r.h.s. is low, and sometimes this factor could yield more accurate results. However, we can neglect the denominator because it gives a higher-order contribution. Finally, the following expression is derived:
\begin{align}
\left[C^{(n)}\right]_{ij} 
&= - \frac{d_{ij}}{d_{ii} d_{jj}}
+ \sum_{\substack{k=1\\k\neq i,j}}^n \frac{d_{ik}d_{kj}}{(1-h\ell_{ii})(1-h\ell_{kk})(1-h\ell_{jj})}
+ \mathcal{O}(h^3) \nonumber \\
&= - \frac{d_{ij}}{d_{ii} d_{jj}}
+ \sum_{\substack{k=1\\k\neq i,j}}^n d_{ik}d_{kj}
+ \mathcal{O}(h^3),
\label{eq_appendix_cij}
\end{align}
which corresponds to \eqref{eq_CN_inverse_2}.


\end{widetext}

\begin{acknowledgments}
This work was supported by JST FOREST Program (Grant Number JPMJFR216K, Japan).
\end{acknowledgments}


\end{document}